\documentclass[12pt]{article}
\usepackage{amssymb}
\usepackage{latexsym,bm}
\usepackage{graphicx}
\usepackage{amsmath}
\usepackage{mathrsfs}
\usepackage{mathrsfs,amscd,amssymb,amsthm,amsmath,bm,graphicx,psfrag,subfigure,url,xcolor}
\usepackage{tikz}

\usepackage{epstopdf}

\date{}
\newcounter{mathitem}
{\begin{list}{{$(\roman{mathitem})$}}{
\setcounter{mathitem}{0}
\usecounter{mathitem}
\setlength{\topsep}{0pt plus 2pt minus 0pt}
\setlength{\parskip}{0pt plus 2pt minus 0pt}
\setlength{\partopsep}{0pt plus 2pt minus 0pt}
\setlength{\parsep}{0pt plus 2pt minus 0pt}
\setlength{\leftmargin}{35pt}
\setlength{\itemsep}{0pt plus 2pt minus 0pt}}}
{\end{list}}
	
\begin{document}
\title{A proof of Bickle's conjecture on collapsible graphs}
\author{\hskip -10mm Xingzhi Zhan}
\maketitle
\footnotetext[1]{Department of Mathematics,  Key Laboratory of MEA (Ministry of Education)
 and Shanghai Key Laboratory of PMMP, East China Normal University, Shanghai 200241, China}
\footnotetext[2]{E-mail address: zhan@math.ecnu.edu.cn}

\begin{abstract}
A graph $G$ is said to be $k$-collapsible if $G$ has minimum degree $k$ and every non-null proper induced subgraph of $G$ has minimum degree less than $k.$
In 2018, Bickle conjectured that the minimum number of vertices of degree $k$ in a $k$-collapsible graph of order $n$ with $k\ge 3$ is
${\rm max}\{\lceil 2n/(2k-1)\rceil,\, k^2-k-2-(k-3)n\}.$ We prove this conjecture.
\end{abstract}

{\bf Keywords.} $k$-collapsible graph; $k$-degenerate graph; minimum degree

{\bf 2020 Mathematics Subject Classification.} 05C07, 05C30, 05C35
\vskip 8mm

\section{Introduction}

We consider finite simple graphs and use standard terminology and notation from [1] and [5]. The {\it order} of a graph is its number of vertices, and the
{\it size} is its number of edges. We denote by $V(G),$ $E(G),$ and $\delta (G)$  the vertex set, edge set, and minimum degree of a graph $G,$ respectively.
The order of $G$ is denoted by $|G|.$  For $S\subseteq V(G)$, let $G[S]$ be the subgraph induced by $S$.  A graph is {\it non-null} if it has at least one vertex.  A graph $G$ is {\it $d$-degenerate} if every non-null induced subgraph of $G$ has a vertex
of degree at most $d$; equivalently, its vertices admit an ordering in which they can be deleted successively with current degree at most $d$. This important concept was introduced by
Lick and White [7] in 1970. We will study a related concept.

{\bf Definition 1.} Let $k$ be a nonnegative integer.  A graph $G$ is {\it $k$-collapsible} if
\[
 \delta(G)=k
 \quad\text{and}\quad
 \delta(F)<k
 \quad\text{for every non-null proper induced subgraph }F\subset G.
\]
$K_1$ is the only $0$-collapsible graph; $K_2$ is the only $1$-collapsible graph; and cycles are the only $2$-collapsible graphs. Any $k$-collapsible graph is necessarily connected.
Bickle introduced the term {\it k-collapsible graph} in his 2010 doctoral dissertation [1, Definition 93],  but the underlying concept had appeared earlier in connection with Palm’s notion of a minimal persistent set at threshold $k$; in later $k$-core terminology, it is a minimal $k$-core.

In 1981, Palm defined a {\it minimal persistent set} at threshold $k$ [8]. In the setting of a simple graph $G$, a vertex set $S\subseteq V(G)$ is persistent at threshold $k$ if
$\delta(G[S])\ge k.$ It is minimally persistent if no proper subset of $S$ is persistent. Minimality forces $\delta(G[S])=k$. Therefore, $G[S]$ is precisely a $k$-collapsible graph. This equivalence was later expressed explicitly in graph-theoretic terminology by Wood and Hicks [10], who called such a graph a {\it minimal $k$-core}. There is also an equivalent formulation in terms of the coloring number [6].

The purpose of this paper is to prove the following result, which was conjectured by Bickle [2, Conjecture 18] in 2018.

{\bf Theorem 1.} {\it The  minimum number of vertices of degree $k$ in a $k$-collapsible graph of order $n$ with $k\ge 3$ is
${\rm max}\{\lceil 2n/(2k-1)\rceil,\, k^2-k-2-(k-3)n\}.$
}

Bickle [2] proved the cases $k=3, 4$ of Theorem 1.

For integers $k\ge3$ and $n\ge k+1$, define
\[
 A(k,n)=\left\lceil\frac{2n}{2k-1}\right\rceil,
 \quad
 B(k,n)=k^2-k-2-(k-3)n,
 \quad
 M(k,n)=\max\{A(k,n),\,B(k,n)\}.
\]

For a graph $G$ of minimum degree $k$, put
\[
 \mathscr{B}(G)=\{v\in V(G)|\,d_G(v)=k\},
 \qquad s(G)=|\mathscr{B}(G)|.
\]
The vertices in $\mathscr{B}(G)$ are called \emph{bottom vertices}; all
other vertices are called \emph{nonbottom vertices}.  The induced
subgraph $G[\mathscr{B}(G)]$ is called the \emph{bottom graph} of $G$.

For a graph $G$, let $d_G(v)$ and $N_G(v)$ denote the degree and neighborhood of $v$, respectively. For vertices $x, y\in V(G),$ we denote by
$\operatorname{dist}_G(x,y)$ the distance between $x$ and $y;$ an {\it $(x, y)$-path} is a path with endpoints $x$ and $y.$

Let $G$ and $H$ be graphs with $V(H)\subseteq V(G)$, and let $u$ be either a vertex or an edge of $G$. We write $u\notin H$ to mean $u\notin V(H)$ if $u$ is a vertex,
and $u\notin E(H)$ if $u$ is an edge.

In Section 2 we prove some auxiliary results, and in Section 3 we prove Theorem 1.

\section{Auxiliary results}

We shall use the following elementary comparison.

{\bf Lemma 2.} {\it If $k\ge5$ and $n\ge k+1$, then
\[
 M(k,n)=
 \begin{cases}
  k+1, & n=k+1,\\
  4,   & n=k+2,\\
  A(k,n), & n\ge k+3.
 \end{cases}
\]
}

{\bf Proof.}
Direct calculation gives
\[
 B(k,k+1)=k+1,
 \qquad
 B(k,k+2)=4.
\]
For $k\ge5$,
\[
 A(k,k+1)=A(k,k+2)=2,
\]
so the first two assertions follow.  Since $k-3>0$, the function
$B(k,n)$ decreases as $n$ increases.  Hence, for $n\ge k+3$,
\[
 B(k,n)\le B(k,k+3)=7-k\le2.
\]
On the other hand, $2n>2k-1$ when $n\ge k+3$, so $A(k,n)\ge2$.
Therefore $A(k,n)\ge B(k,n)$ throughout this range.
\hfill $\Box$

{\bf Lemma 3.} {\it Let $G$ be a graph with $\delta(G)=k$ where $k\ge 1.$  Then $G$ is
$k$-collapsible if and only if $G-v$ is $(k-1)$-degenerate for every $v\in V(G)$.
}

{\bf Proof.}
Suppose first that $G$ is $k$-collapsible.  Every non-null induced
subgraph of $G-v$ is a proper induced subgraph of $G$ and hence has
minimum degree at most $k-1$.  Thus $G-v$ is $(k-1)$-degenerate.

Conversely, suppose that $G-v$ is $(k-1)$-degenerate for every
$v\in V(G)$.  If a proper induced subgraph $F$ of $G$ had minimum
degree at least $k$, then choosing $v\in V(G)\setminus V(F)$ would make
$F$ an induced subgraph of $G-v$, contradicting the
$(k-1)$-degeneracy of $G-v$.
\hfill $\Box$

The next criterion will certify the first construction family.

{\bf Lemma 4.} {\it Let $G$ be a graph, let $k$ be a positive integer, and suppose that $X,Y\subseteq V(G)$ partition $V(G)$.
Assume that
\begin{enumerate}
 \item[(1)] $G[X]$ is connected and $d_G(x)=k$ for every $x\in X;$
 \item[(2)] $ d_G(y)\ge k+1$ and $N_G(y)\cap X\ne\varnothing$ for every $y\in Y;$
 \item[(3)] $G[Y]$ is $(k-1)$-degenerate.
\end{enumerate}
Then $G$ is $k$-collapsible and $X=\mathscr{B}(G)$.}

{\bf Proof.}
Conditions (1) and (2) imply that $\delta(G)=k$ and
$X=\mathscr{B}(G)$.  Fix $v\in V(G)$.  We exhibit a deletion ordering
of $G-v$ in which every vertex has current degree at most $k-1$.

First suppose that $v\in X$.  Let $T$ be a spanning tree of $G[X]$ and
write
\[
 X\setminus\{v\}=\{x_1,x_2,\ldots,x_q\},
\]
where
\[
 \operatorname{dist}_T(v,x_1)\le\operatorname{dist}_T(v,x_2)
 \le\cdots\le\operatorname{dist}_T(v,x_q).
\]
For every $i$, let $y_i$ be the neighbor of $x_i$ on the unique
$(v, x_i)$-path in $T$ that is closer to $v$.  Then $y_i=v$, or
$y_i=x_j$ for some $j<i$.  Thus $y_i\notin (G-v)-\{x_1,x_2,\ldots,x_{i-1}\}.$
Since $d_G(x_i)=k$, the current degree of $x_i$ in $(G-v)-\{x_1,x_2,\dots,x_{i-1}\}$
is at most $k-1$. Hence all vertices of $X\setminus\{v\}$ can be deleted in the displayed
order.

Now suppose that $v\in Y$.  By (2), choose $r\in N_G(v)\cap X$.
Since $d_G(r)=k$ and $rv\in E(G)$, $d_{G-v}(r)=k-1,$ so delete $r$ first.  Let $S$ be a spanning tree of $G[X]$ and write
\[
 X\setminus\{r\}=\{w_1,w_2,\ldots,w_q\},
\]
where
\[
 \operatorname{dist}_S(r,w_1)\le\operatorname{dist}_S(r,w_2)
 \le\cdots\le\operatorname{dist}_S(r,w_q).
\]
For each $i$, let $u_i$ be the neighbor of $w_i$ on the unique
$(r, w_i)$-path in $S$ that is closer to $r$.  Then $u_i=r$, or
$u_i=w_j$ for some $j<i$.  Consequently
\[
 u_i\notin (G-v)-r-\{w_1,w_2,\ldots,w_{i-1}\},
\]
and the current degree of $w_i$ in $(G-v)-r-\{w_1,w_2,\ldots,w_{i-1}\}$
is at most $d_G(w_i)-1=k-1$.  Thus every vertex of $X$ can be deleted.

After the vertices of $X\cap V(G-v)$ have been deleted, the remaining
graph is
\[
 G[Y]\quad\text{if }v\in X,
 \qquad\text{or}\qquad
 G[Y\setminus\{v\}]\quad\text{if }v\in Y.
\]
Both graphs are $(k-1)$-degenerate by (3), because an induced subgraph
of a $(k-1)$-degenerate graph is again $(k-1)$-degenerate.  Therefore
all remaining vertices can be deleted with current degree at most
$k-1$.  Lemma 3 completes the proof. \hfill $\Box$

A finite sequence of nonnegative integers $\mathbf d=(d_1,d_2,\ldots,d_n)$
is called a \emph{graphic sequence}  if there exists a simple graph $G$ with vertices
$v_1,v_2,\ldots,v_n$ such that
\[
  \deg_G(v_i)=d_i \qquad (1\le i\le n).
\]
Such a graph $G$ is called a \emph{realization} of $\mathbf d.$
A \emph{connected realization} of $\mathbf d$ is a realization of
$\mathbf d$ that is connected.

{\bf Lemma 5.} {\it Let $k\ge5$ and $h\ge k+1$.  Every even-sum sequence of length $h$
whose entries belong to $\{k-1,k\}$ has a connected realization.}

{\bf Proof.} Write the nonincreasing sequence as
\[
 \bigl(\underbrace{k,\ldots,k}_{a\text{ times}},
       \underbrace{k-1,\ldots,k-1}_{h-a\text{ times}}\bigr).
\]
The Erd\H{o}s--Gallai criterion ([4] or [3, p.11])  states that a nonincreasing
sequence $d_1,\ldots,d_h$ of nonnegative integers is graphic if and
only if the sum of its terms is even and
\[
 \sum_{i=1}^{r}d_i
 \le r(r-1)+\sum_{i=r+1}^{h}\min\{r,d_i\}
 \qquad\text{for every }r\in\{1,2,\ldots,h\}.       \eqno (1)
\]
The sequence under consideration has even sum by hypothesis.  We now
verify the inequalities in (1).
Define the Erd\H{o}s--Gallai difference at $r$ by
$$
F(r)=r(r-1)+\sum_{i=r+1}^{h}\min\{r,d_i\}-\sum_{i=1}^{r}d_i.
$$
The inequality (1) is equivalent to $F(r)\ge0$.
For an index $r\le k-1$, the right-hand side is $r(h-1)$, while the left-hand side is at most
$rk\le r(h-1)$. Hence $F(r)\ge 0.$

For $r\ge k$, let $S$ be the sum of the terms of the sequence, let
$P_r$ be the sum of the first $r$ entries, and then we have
\[
 F(r)=r(r-1)+S-2P_r.
\]
Moreover,
\[
 F(r+1)-F(r)=2r-2d_{r+1}\ge0,
\]
so it is enough to check $r=k$.  Since
\[
 P_k=k(k-1)+\min\{a,k\},
 \qquad
 S=(k-1)h+a,
\]
we have
\[
 F(k)=(k-1)(h-k)+a-2\min\{a,k\}.
\]
If $a<k$, then $F(k)\ge k-1-a\ge0$.  If $a\ge k$ and
$h\ge k+2$, then $F(k)\ge k-2\ge0$.  Finally, if $h=k+1$, then
$a\in\{k,k+1\}$ in the remaining cases; $a=k+1$ gives equality,
whereas for $a=k$ the sum of the terms is odd, contrary to the
hypothesis.  Thus the sequence is graphic.

Every realization of this sequence has minimum degree at least $k-1\ge4$.  Suppose that
a realization has components $C_1,\ldots,C_r$, where $r\ge2$.  Every
$C_i$ contains a cycle and hence contains an edge that is not a
cut-edge.  Choose such edges $x_1y_1\in E(C_1)$ and
$x_2y_2\in E(C_2)$.  Delete $x_1y_1$ and $x_2y_2$, and add
$x_1x_2$ and $y_1y_2$.  Since the deleted edges are not cut-edges,
both $C_1-x_1y_1$ and $C_2-x_2y_2$ remain connected.  The two new
edges therefore merge $C_1$ and $C_2$ into one component.  Moreover,
the resulting graph is simple, because there were no edges between
$C_1$ and $C_2$, and the operation does not change the degree of any
vertex.  Thus the number of components decreases from $r$ to $r-1$,
while the degree sequence is preserved.  Regard the merged component
as a new component; it still has minimum degree at least $k-1\ge4$.
Repeating this operation reduces the number of components by exactly
one at each step, so after exactly $r-1$ steps we obtain a connected
realization.
\hfill $\Box$

For integers $n>k\ge2$ with $kn$ even, denote by $H_{k,n}$ the Harary
graph with vertex set $\{1,2,\ldots,n\}$, with cyclic distance computed modulo $n$.
If $k=2q$, join vertices at cyclic distance at most $q$.  If
$k=2q+1$, then $n$ is even; join vertices at cyclic distance at most
$q$ and also join antipodal pairs.  The graph $H_{k,n}$ is $k$-regular
and $k$-connected ([5] or [9, p.150]).  We write it as $H_{k,h}$ when its
order is denoted by $h$.

{\bf Lemma 6.} {\it Let $k\ge5$ and $h\ge k+1$ be integers, and suppose that $kh$ is even.  Suppose that either
(1) $k$ is even and $h\ge2k-3$, or (2) $k$ is odd and $h\ge3k-5$.
Then there is a graph $Q$ of order $h$ with specified vertices $c,z$
such that
\[
 d_Q(c)=k-1,
 \qquad d_Q(z)=k+1,
\]
all other vertices have degree $k$, and $Q-z$ is connected.}

{\bf Proof.}
Start with $R=H_{k,h}$ and label $c=1$, $u=2$, and $z=k+1$.  The edge
$cu$ is present.  If $k$ is even, the two directed cyclic distances
from $u$ to $z$ are $k-1$ and $h-k+1$.  The latter is at least $k-2$,
and both numbers exceed $k/2$.  If $k$ is odd, then the cyclic distance
between $u$ and $z$ is $k-1$, which is greater than $(k-1)/2$ and less
than $h/2$.  Consequently, $u$ and $z$ are not antipodal.  Thus $uz\notin R$
in both cases.  Define a graph $Q$ by
\[
 V(Q)=V(R),
 \qquad
 E(Q)=\bigl(E(R)\setminus\{cu\}\bigr)\cup\{uz\}.
\]
The stated degrees follow.  Since $R$ is $k$-connected, $R-z$ is
$(k-1)$-connected and therefore has edge connectivity at least
$k-1\ge4$.  Deleting the single edge $cu$ cannot disconnect $R-z$,
so $Q-z$ is connected.
\hfill $\Box$

For each fixed integer $k\ge5$, define $\mathbf a_h$ and
$\mathbf b_h$ for every integer $h\ge k+1$, and define $\mathbf c_h$
for every integer $h\ge k+2$, as follows:
\begin{align*}
 \mathbf a_h
   &=\bigl(\underbrace{k,\ldots,k}_{h-4\text{ times}},
      k-1,k-1,k-1,k-2\bigr),\\
 \mathbf b_h
   &=\bigl(\underbrace{k,\ldots,k}_{h-5\text{ times}},
      k-1,k-1,k-1,k-1,k-2\bigr),\\
 \mathbf c_h
   &=\bigl(k+1,
      \underbrace{k,\ldots,k}_{h-5\text{ times}},
      k-1,k-1,k-1,k-2\bigr).
\end{align*}

We first construct several auxiliary graphs, which we call units.

{\bf Lemma 7.} {\it Let $k\ge5$ be an integer.
\begin{enumerate}
 \item[(1)] If $k$ is odd, then $\mathbf a_{2k-3}$ is graphic, and every
       realization of $\mathbf a_{2k-3}$ is connected.  For every $k\ge 5$,
       the sequence $\mathbf b_{2k-4}$ is graphic, and every realization
       of $\mathbf b_{2k-4}$ is connected.
 \item[(2)] If $k$ is even, then $\mathbf c_{2k-3}$ has a realization $Q$
       such that $Q-z$ is connected, where $z$ is its unique
       degree-$(k+1)$ vertex.  For every $k\ge 5$, the sequence
       $\mathbf c_{3k-5}$ has such a realization.
\end{enumerate}
}

{\bf Proof.} For (1), the sums of the two sequences are respectively
\[
 k(2k-3)-5
 \qquad\text{and}\qquad
 k(2k-4)-6,
\]
which are even under the stated hypotheses.  Consider a sequence
\[
 D=\bigl(\underbrace{k,\ldots,k}_{a\text{ times}},
          \underbrace{k-1,\ldots,k-1}_{b\text{ times}},k-2\bigr)
\]
of length $h=a+b+1$.  Let $S$ be the sum of its terms and $P_r$ the
sum of its first $r$ entries.  We apply the Erd\H{o}s--Gallai
criterion.  For $r\le k-2$, the right-hand side is $r(h-1)$ and the left-hand side is $P_r\le rk$.  For $r\ge k$, the Erd\H{o}s--Gallai difference
\[
 F(r)=r(r-1)+S-2P_r
\]
is nondecreasing.  It therefore remains only to check $r=k-1$ and
$r=k$.  Direct substitution gives the two differences
\begin{align*}
 F(k-1)&=(k-1)(h-1)-1-(k-1)^2-\min\{a,k-1\},\\
 F(k)&=S-k(k-1)-2\min\{a,k\}.
\end{align*}

For $\mathbf a_{2k-3}$, substitute
\[
 h=2k-3,
 \qquad a=2k-7,
 \qquad S=k(2k-3)-5.
\]
When $k=5$, both differences equal $4$; for odd $k\ge7$, they are
\[
 (k-1)(k-4)-1
 \qquad\text{and}\qquad
 k^2-4k-5,
\]
respectively.  For $\mathbf b_{2k-4}$, substitute
\[
 h=2k-4,
 \qquad a=2k-9,
 \qquad S=k(2k-4)-6.
\]
For $k=5,6,7,8$, the pairs of differences are respectively
\[
 (2,2),\ (6,6),\ (12,12),\ (20,20),
\]
and for $k\ge9$ they are
\[
 (k-1)(k-5)-1
 \qquad\text{and}\qquad
 k^2-5k-6.
\]
All these differences are nonnegative, so both sequences are graphic.

Their minimum degree is $k-2$.  Every component of a realization
therefore has at least $k-1$ vertices, while the total order is less
than $2k-2$.  Hence every realization is connected.

For (2), let $h=2k-3$ or $h=3k-5$, as applicable, and start with $R=H_{k,h}$.
The vertices labeled below are distinct, and $z=k+3$ exists for both of these orders.
Label
\[
 x=1,
 \qquad c=2,
 \qquad u=3,
 \qquad v=4,
 \qquad w=5,
 \qquad z=k+3.
\]
The pairs $xc,xu,vw$ are edges of $R$.  If $k$ is even, then
$h\in\{2k-3,3k-5\}$, and the two directed cyclic distances between
$z$ and $w$ are $k-2$ and $h-k+2$, both greater than $k/2$.  If $k$
is odd, then $h=3k-5$, and the cyclic distance between $z$ and $w$ is
$k-2$, which is greater than $(k-1)/2$ and less than $h/2$.
Consequently, $z$ and $w$ are not antipodal.  Thus $zw$ is a nonedge.
Define a graph $Q$ by
\[
 V(Q)=V(R),
 \qquad
 E(Q)=\bigl(E(R)\setminus\{xc,xu,vw\}\bigr)\cup\{zw\}.
\]
The resulting degrees are $k-2$ at $x$, $k-1$ at $c,u,v$, $k+1$ at
$z$, and $k$ elsewhere.  Since $R-z$ is $(k-1)$-connected, its edge
connectivity is at least $k-1\ge4$.  Deleting the three displayed
edges does not disconnect it, so $Q-z$ is connected.
\hfill $\Box$

\section{Proof of the main theorem}

Theorem 1 states that for all integers $k\ge 3$ and $n\ge k+1$,
\[
 \min\{s(G):|G|=n\text{ and }G\text{ is }k\text{-collapsible}\}
 =M(k,n). \eqno (2)
\]

{\bf Proof of Theorem 1.} Denote by $f(k, n)$ the number on the left-hand side of (2).

Bickle has proved that
$$
f(k,n)\ge M(k,n) \eqno (3)
$$
[2, Theorem 14 and Lemma 15] and that for $k=3, 4,$ equality holds in (3) [2, Theorems 16 and 17].
Thus, to prove (2) it suffices to show that for all pairs of integers $(k, n)$ with $k\ge 5$ and $n\ge k+1,$
there exists a $k$-collapsible graph $G$ with exactly $M(k,n)$ bottom vertices.

If $n=k+1$, take $G=K_{k+1}$.  Every vertex has degree $k$, and every proper induced subgraph has order at
most $k$ and hence minimum degree at most $k-1$.  Thus $K_{k+1}$ is $k$-collapsible and has $k+1=M(k,k+1)$ bottom vertices
by Lemma 2.

Suppose next that $n=k+2$.  Define $J_{k+2}$ explicitly by
\[
 V(J_{k+2})=\{1,2,\ldots,k+2\}
\]
and
\[
 E(J_{k+2})
 =\{(i,j):1\le i<j\le k+2\}\setminus\{(1,2),(3,4)\}.
\]
Thus $J_{k+2}$ is the graph sometimes denoted by $K_{k+2}-2K_2$.
The four endpoints of the two deleted edges have degree $k$, and all
other vertices have degree $k+1$.  A proper induced subgraph with
minimum degree at least $k$ would have exactly $k+1$ vertices and would
have to be complete.  However, deleting one vertex from $J_{k+2}$
cannot remove both missing independent edges.  Therefore $J_{k+2}$ is
$k$-collapsible and has four bottom vertices.  By Lemma 2, $4=M(k,k+2)$.

It remains only to consider $n\ge k+3\ge8.$ By Lemma 2,
\[
 s:=M(k,n)=A(k,n)=\left\lceil\frac{2n}{2k-1}\right\rceil
\]
is the required number of bottom vertices. We distinguish two cases.

\textbf{Case 1. $n\le(k-1)s+1$ (the connected-bottom construction)}

Put $h=n-s$. The case assumption is
\[
 n\le(k-1)s+1.                                      \eqno (4)
\]
Put $t=(k-1)s+2-n.$ We first establish the parameter bounds
\[
 h\ge k+1,
 \qquad
 1\le t\le h-2.                                     \eqno (5)
\]
To prove $h\ge k+1$, it is enough to show that
\[
n-k-1\ge \frac{2n}{2k-1}.
\]
Indeed, the left-hand side is an integer and
$s=\lceil2n/(2k-1)\rceil$.  Moreover,
\begin{align*}
 (2k-1)(n-k-1)-2n
   &=(2k-3)n-(2k-1)(k+1)\\
   &\ge(2k-3)(k+3)-(2k-1)(k+1)\\
   &=2k-8\ge0.
\end{align*}
Consequently
\[
 n-k-1\ge \left\lceil\frac{2n}{2k-1}\right\rceil=s,
\]
and hence $h=n-s\ge k+1$.

The inequality (4) immediately gives
\[
 t=(k-1)s+2-n\ge1.
\]
It remains to prove that $t\le h-2;$ i.e., $h-t\ge2$.  Observe that $h-t=2n-ks-2.$
We have $s\ge2$ because $n\ge k+3$.  If $s=2$, then
\[
 h-t=2n-2k-2\ge2(k+3)-2k-2=4\ge 2.
\]
Now suppose that $s\ge3$.  The definition of $s$ gives $(s-1)(2k-1)<2n,$ where both sides are integers, and therefore
\[
 2n\ge(s-1)(2k-1)+1.
\]
It follows that
\begin{align*}
 h-t
   &=2n-ks-2\\
   &\ge(s-1)(2k-1)+1-ks-2\\
   &=(k-1)s-2k\\
   &\ge3(k-1)-2k\\
   &=k-3\ge2.
\end{align*}
This proves (5).

We now specify two vertex sets and all edges between them explicitly.
Let
\[
 X=\{x_1,x_2,\ldots,x_s\},
 \qquad
 Y=\{y_1,y_2,\ldots,y_h\},
\]
and let the edges within $X$ be $E_X=\{(x_i,x_{i+1}):1\le i\le s-1\}.$ Thus the graph on $X$ with edge set $E_X$ is the path
$x_1x_2\cdots x_s$.  Define
\[
 r_i=
 \begin{cases}
  k-1, & i\in\{1,s\},\\
  k-2, & 2\le i\le s-1,
 \end{cases}
 \qquad
 q_0=0,
 \qquad
 q_i=\sum_{a=1}^{i}r_a.
\]
Then
\[
 q_s=2(k-1)+(s-2)(k-2)=(k-2)s+2=h+t.
\]
For every positive integer $q$, put
\[
 \rho(q)=1+((q-1)\bmod h),
\]
and define the edges between $X$ and $Y$ by
\[
 E_{XY}=\{(x_i,y_{\rho(q)}):1\le i\le s,
                    \ q_{i-1}<q\le q_i\}.          \eqno (6)
\]

This gives the required incidence pattern.  Indeed, the interval
$q_{i-1}<q\le q_i$ has length $r_i$, and
\[
 r_i\le k-1\le h-2<h.
\]
Consequently the values $\rho(q)$ in this interval are distinct, so
$x_i$ has exactly $r_i$ distinct neighbors in $Y$. Therefore (6) defines a simple bipartite graph.
Thus each endpoint of the path on $X$ has $k-1$ neighbors in $Y$, and each internal vertex has
$k-2$ such neighbors.

Moreover, $q$ runs through $1,2,\ldots,h+t$, where $1\le t\le h-2$.
Hence the residues $1,\ldots,t$ occur twice, whereas the residues
$t+1,\ldots,h$ occur once.  Thus $y_1,\ldots,y_t$ each receive two
incident edges in (6), whereas $y_{t+1},\ldots,y_h$ each receive one.
 Call
$y_1,\ldots,y_t$ \emph{double} vertices and
$y_{t+1},\ldots,y_h$ \emph{single} vertices.  By (5), there is at
least one double vertex and at least two single vertices.

We next define the graph with vertex set $Y$.  Initially prescribe internal
degree $k-1$ for every double vertex and internal degree $k$ for every
single vertex.  The resulting degree sequence is
\[
 \bigl(\underbrace{k,\ldots,k}_{h-t\text{ times}},
       \underbrace{k-1,\ldots,k-1}_{t\text{ times}}\bigr).       \eqno (7)
\]
If its sum is even, let $Q$ be a connected realization on $Y$ supplied
by Lemma 5.  If its sum is odd and $t\ge2$, raise the prescribed internal degree of $y_1$ from $k-1$ to $k$ and
let $Q$ be a connected realization of the new even-sum sequence.  The double vertex $y_2$ still has internal degree $k-1$.

It remains to treat the case in which the sum in (7) is odd and $t=1$.
The sum in (7) is $kh-t=kh-1$, so $kh$ is even.  Use Lemma 6 to
obtain a graph $Q$ on $Y$, relabeled so that its degree-$(k-1)$ vertex is $c=y_1$,
the unique double vertex, and its degree-$(k+1)$ vertex is $z=y_2$, a single vertex.  The hypotheses of Lemma 6 hold.
Indeed, $t=1$ gives
\[
 n=(k-1)s+1,
 \qquad
 h=(k-2)s+1.
\]
If $k$ is even, $s\ge2$ gives $h\ge2k-3$.  If $k$ is odd, then $kh$
even implies that $h$ is even.  Since $h=(k-2)s+1$, this forces $s$
to be odd.  Hence $s\ge3$ and $h\ge3k-5$.

In every case, define the required graph $G$ explicitly by
\[
 V(G)=X\cup Y,
 \qquad
 E(G)=E_X\cup E_{XY}\cup E(Q).
\]

We now verify that $G$ has the required properties.  Every
vertex of $X$ has degree $k$.  An unadjusted double vertex has total
degree $(k-1)+2=k+1$, every single vertex has total degree $k+1$, and
the adjusted double vertex, when present, has total degree $k+2$.  In
the exceptional construction, $c$ has total degree $k+1$, $z$ has
total degree $k+2$, and every other vertex of $Y$ has total degree
$k+1$.  Thus $X$ is exactly the bottom set, $G[X]=P_s$ is connected,
and every vertex of $Y$ has a neighbor in $X$.

It remains to show that $G[Y]$ is $(k-1)$-degenerate by Lemma 4.
In each of the two nonexceptional cases---that is, when the sum in (7)
is even, or when it is odd and $t\ge 2$---the graph $G[Y]$ is connected
and has maximum degree at most $k$.  If the sequence in (7) has even sum, choose any double
vertex as $c$; if its sum is odd and $t\ge2$, let $c=y_2$.  In either
case $d_{G[Y]}(c)=k-1$.  Choose a spanning tree
rooted at $c$.  Delete $c$ first, and then delete the other vertices
in nondecreasing order of their distance from $c$ in the tree.  Every
vertex after $c$ has lost the edge to the preceding vertex on its tree
path to $c$, so its current degree is at most $k-1$.

In the exceptional case, apply the same distance ordering in the
connected graph $G[Y]-z$, starting at $c$.  Every vertex of this graph
has degree at most $k$, and $c$ has degree at most $k-1$.  After all vertices
of $G[Y]-z$ have been deleted, $z$ is isolated and can be deleted.
Thus $G[Y]$ is $(k-1)$-degenerate in every case.
Lemma 4 now shows that $G$ is $k$-collapsible.  Since $X=\mathscr{B}(G)$ and $|X|=s=M(k,n)$, this
graph has the required number of bottom vertices.

\textbf{Case 2. $n>(k-1)s+1$ (the cyclic-unit construction)}

The case assumption is
\[
 n>(k-1)s+1.                                       \eqno (8)
\]
Put
\[
 N=\left\lfloor\frac{s(2k-1)}{2}\right\rfloor,
 \qquad
 p=\left\lfloor\frac{s}{2}\right\rfloor,
 \qquad
 d=N-n.
\]
Write $n=(k-1)s+j$.  The following identities will be used in the
construction.  Since $s=\lceil2n/(2k-1)\rceil$, we have
\[
 n\le\left\lfloor\frac{s(2k-1)}{2}\right\rfloor=N.
\]
Moreover,
\[
 N=(k-1)s+\lfloor s/2 \rfloor=(k-1)s+p.
\]
Condition (8) gives $j\ge2$, while $n\le N$ gives $j\le p$.
Consequently $d=N-n=p-j$, and hence $0\le d\le p-2$.  In particular,
$p\ge2$.

In the construction below, the vertices of each unit are designated
as bottom or nonbottom vertices.  The \emph{designated bottom graph}
of a unit is the subgraph induced by its designated bottom vertices.

We use three types of units.  Each unit contains a connected designated bottom
graph and a connected internal graph $Q$ on its nonbottom vertices.  The edges
between these two parts are exactly the edges specified below.

\medskip
\noindent
\emph{Full pair unit.}
The bottom graph is the edge $ab$, and $Q$ has $2k-3$ vertices.  Choose
$c\in V(Q)$ and join it to both $a$ and $b$.  Partition
$V(Q)\setminus\{c\}$ into two sets of size $k-2$ and join their
vertices respectively to $a$ and $b$.  The unit has order $2k-1$.

\medskip
\noindent
\emph{Short pair unit.}
The bottom graph is the edge $ab$, and $Q$ has $2k-4$ vertices.  Choose
distinct vertices $c,c'\in V(Q)$ and join both of them to both $a$ and
$b$.  Partition the remaining vertices into two sets of size $k-3$
and join their vertices respectively to $a$ and $b$.  The unit has
order $2k-2$.

\medskip
\noindent
\emph{Full triple unit.}
The bottom graph is the path $a_1a_2a_3$, and $Q$ has $3k-5$ vertices.
Choose $c\in V(Q)$ and join it to both $a_1$ and $a_3$.  Partition
$V(Q)\setminus\{c\}$ into three sets of size $k-2$ and join their
vertices respectively to $a_1,a_2,a_3$.  The unit has order $3k-2$.

\medskip
The graph $Q$ in a full pair unit is a connected realization of  $\mathbf a_{2k-3}$ when $k$ is odd and
of $\mathbf c_{2k-3}$ when $k$ is even.  The graph $Q$ in a short pair unit is a connected realization of
$\mathbf b_{2k-4}$, and the graph $Q$ in a full triple unit is a connected realization of $\mathbf c_{3k-5}$.
The existence of these graphs $Q$ is guaranteed by Lemma 7.

For each sequence, label the unique degree-$(k-2)$ vertex $x$ and
three degree-$(k-1)$ vertices $c,u,v$.  In a short pair unit, label the
fourth degree-$(k-1)$ vertex $c'$.  In a $\mathbf c$-unit, label the
degree-$(k+1)$ vertex $z$.  Assign $x,u,v,z$, whenever present, to
single-attachment positions in the corresponding unit, and assign
$c,c'$ to the double-attachment positions.

There are enough single-attachment positions for these special
vertices.  A full pair unit has $2k-4$ such positions, a short pair
unit has $2k-6$, and a full triple unit has $3k-6$.

For each unit, choose its internal graph $Q$ from Lemma 7, as prescribed above, and relabel its vertices
to match the specified attachment roles.  This completes the
definitions of the three unit types.

If $s=2p$, begin with $p$ full pair units.  If $s=2p+1$, begin with
$p-1$ full pair units and one full triple unit.  In either case the
total order is $N=\lfloor s(2k-1)/2\rfloor.$ By the preceding calculation, replace exactly $d$ full pair units by
short pair units.  There are enough full pair units because
$d\le p-2$.  Each replacement reduces the order by one.  The resulting
$p$ units therefore have total order $n$ and exactly $s$ designated
bottom vertices.

Take pairwise vertex-disjoint copies of these $p$ units and index them
cyclically, with subscripts read modulo $p$:
\[
 U_1,U_2,\ldots,U_p.
\]
In $U_i$, attach the subscript $i$ to the labels introduced above; in
particular, its relevant nonbottom vertices are denoted by
$x_i,u_i,v_i$, and, when present, $c_i,c'_i,z_i$.  Let
\[
 E_{\times}
 =\{(u_i,x_{i+1}),(v_i,x_{i+1}):1\le i\le p\},   \eqno (9)
\]
where $x_{p+1}=x_1$.  The target graph $G$ is obtained from the
disjoint union $U_1+U_2+\cdots+U_p$ by adding the edges in
$E_{\times}$.  Equivalently,
\[
 V(G)=\bigcup_{i=1}^{p}V(U_i),
 \qquad
 E(G)=\left(\bigcup_{i=1}^{p}E(U_i)\right)\cup E_{\times}.
\]

We verify the degrees directly.  After the cross edges (9) are added,
every bottom vertex has degree $k$ and every nonbottom vertex has
degree greater than $k$.  In fact, all nonbottom vertices have total degree $k+1$, except that $z_i$,
 when present, has total degree $k+2$.

We next establish the deletion mechanism.  Suppose that one vertex
$w$ of a unit $U_i$ is absent.  We show that all remaining vertices
of $U_i$ can be deleted successively with current degree at most
$k-1$ and that, afterward, $x_{i+1}$ can also be deleted if it is
present.

First delete the remaining bottom vertices of $U_i$.  If $w$ is a
bottom vertex, choose a spanning tree of the original connected
bottom graph rooted at $w$ and order the other bottom vertices by
nondecreasing distance from $w$.  The parent of each vertex in this
order is absent or has already been deleted.  Since every bottom
vertex has degree $k$, its current degree is therefore at most $k-1$.
If $w$ is nonbottom, choose a bottom neighbor $b$ of $w$.  The vertex
$b$ has lost the edge $bw$ and hence has current degree at most
$k-1$.  Delete $b$, choose a spanning tree of the bottom graph rooted
at $b$, and then delete the other bottom vertices in nondecreasing
order of their distance from $b$.  Again, each vertex has an absent
or previously deleted parent and hence has current degree at most
$k-1$ when it is deleted.

After the bottom graph has disappeared, consider the internal graph
$Q_i$.  For a unit based on $\mathbf a_h$ or $\mathbf b_h$, choose a
spanning tree $T$ of the original connected graph $Q_i$, rooted at
$c_i$.  Process the vertices in nondecreasing order of their distance
from $c_i$ in $T$, omitting any vertex that is already absent.  Thus,
if $c_i=w$, the root is simply omitted.  Every processed vertex other
than a present root has an absent or previously deleted parent in
$T$.  When the still-present cross edges are counted, every remaining vertex of
$Q_i$ has current degree at most $k$, while a present $c_i$ has
current degree $k-1$.  Hence every processed vertex has current degree
at most $k-1$.

For a unit based on $\mathbf c_h$, use a spanning tree of the original
connected graph $Q_i-z_i$, rooted at $c_i$, and apply the same ordering,
again omitting any vertex already absent.  Every remaining non-$z_i$ vertex has
current degree at most $k$ before losing the edge to its tree parent,
even while $z_i$ and the cross edges remain present.  After all
vertices of $Q_i-z_i$ have been deleted, $z_i$ is isolated and can be
deleted.  Thus the whole unit disappears in every case.

Finally, deleting $u_i$ and $v_i$ removes the two incoming cross edges
at $x_{i+1}$.  The latter then has degree at most
\[
 (k-2)+1=k-1,
\]
where the remaining term $1$ is its bottom attachment.  Thus
$x_{i+1}$ can be deleted if it is present. When $x_{i+1}$ is deleted, the bottom graph of $U_{i+1}$ has a vertex
$b_{i+1}$ which is a neighbor of $x_{i+1}$ and the degree of $b_{i+1}$ decreases from $k$ to $k-1.$ With $b_{i+1}$
in place of the vertex $w$ above, keep deleting vertices until all vertices of $U_{i+1}$ are deleted.
Cyclically relabel the units, if necessary, so that the unit containing the initially
absent vertex is $U_1$. Perform this deletion process successively for $i=1,2,\dots,p.$

Clearly we have $|G|=n$, and the degree verification identifies the $s$ designated bottom vertices as exactly
the bottom vertices of $G$.  It remains to prove that $G$ is
$k$-collapsible.  Delete an arbitrary vertex $v$.  The deletion mechanism removes the
unit containing $v$ and then triggers the next unit.  Repeating it
around the cyclic order deletes every vertex of $G-v$, always with
current degree at most $k-1$.  Hence $G-v$ is $(k-1)$-degenerate for
every $v\in V(G)$.  Lemma 3 shows that $G$ is $k$-collapsible.  Finally, $s=M(k,n)$, so $G$ has the required number
of bottom vertices. \hfill $\Box$

\section*{\normalsize Declaration of AI Use}

ChatGPT was used to assist in developing and checking the constructions and proofs. The author independently verified all mathematical arguments, wrote the paper, and takes full responsibility for its content.

\vskip 5mm
{\bf Acknowledgements.} This research  was supported by the NSFC Grant No. 12271170 and by Science and Technology Commission of Shanghai Municipality
 Grant No. 22DZ2229014.

\end{document}